\documentclass{amsart}%
\usepackage{amscd}
\usepackage{amsmath}
\usepackage{amsfonts}
\usepackage{amssymb}
\usepackage{graphicx}%
\providecommand{\U}[1]{\protect\rule{.1in}{.1in}}
\newtheorem{theorem}{Theorem}
\theoremstyle{plain}
\newtheorem{acknowledgement}{Acknowledgement}

\newtheorem{corollary}{Corollary}

\newtheorem{definition}{Definition}

\newtheorem{lemma}{Lemma}

\newtheorem{proposition}{Proposition}

\numberwithin{equation}{section}

\begin{document}
\title[Lower Bounds for Maximal operators ]{Local lower norm estimates for dyadic maximal operators and related Bellman fuctions}
\author{Antonios D. Melas}
\author{Eleftherios Nikolidakis}
\address{Department of Mathematics, University of Athens, Panepistimiopolis 15784,
Athens, Greece}
\email{amelas@math.uoa.gr, lefteris@math.uoc.gr}
\date{June 25, 2015}
\subjclass[2010]{ 42B25}
\keywords{Bellman, dyadic, maximal}

\begin{abstract}
We provide lower $L^{q}$ and weak $L^{q}$-bounds for the localized dyadic
maximal operator on $\mathbb{R}^{n}$, when the local $L^{1}$ and the local
$L^{p}$ norm of the function are given. We actually do that in the more
general context of homogeneous tree-like families in probability spaces.

\begin{acknowledgement}
This research has been co-financed by the European Union and Greek national
funds through the Operational Program "Education and Lifelong Learning" of the
National Strategic Reference Framework (NSRF). ARISTEIA\ I, MAXBELLMAN 2760,
research number 70/3/11913.

\end{acknowledgement}
\end{abstract}
\maketitle

\section{Introduction}

The dyadic maximal operator on $\mathbb{R}^{n}$ is a useful tool in analysis
and is defined by
\begin{equation}
M{}\,_{d}\phi{}(x)=\sup\left\{  \frac{1}{\left\vert Q\right\vert }\int
_{Q}\left\vert \phi(u)\right\vert du:x\in Q\text{, }Q\subseteq\mathbb{R}%
^{n}\text{ is a dyadic cube}\right\}  \label{i1}%
\end{equation}
for every $\phi\in L_{\text{loc}}^{1}(\mathbb{R}^{n})$ where the dyadic cubes
are the cubes formed by the grids $2^{-N}\mathbb{Z}^{n}$ for $N=0,1,2,...$.

Localizing the above operator on a unit cube $Q_{0}$ (that is considering only
sucubes of $Q_{0}$ in the above supremum) leads to an operator that can be
generalized in the context of a $(X,\mu)$ be a nonatomic probability space
$(X,\mu)$ equipped with a tree like family (see also \cite{Mel1a}). The
precise definition follows:

\begin{definition}
(a) A set $\mathcal{T}$ of measurable subsets of $X$ will be called an
$N$-homogeneous tree-like family (where $N>1$ is an integer) if the following
conditions are satisfied:

(i) $X\in\mathcal{T}$ \ and for every $I\in\mathcal{T}$ \ there corresponds a
finite subset $\mathcal{C}(I)\subseteq\mathcal{T}$ \ containing $N$ elements
each having measure equal to $N^{-1}\mu(I)$ such that the elements of
$\mathcal{C}(I)$ are pairwise disjoint subsets of $I$ and $I=\bigcup
\mathcal{C}(I)$.

(ii) $\mathcal{T}=\bigcup_{m\geq0}\mathcal{T}_{(m)}$ where $\mathcal{T}%
_{(0)}=\{X\}$ and $\mathcal{T}_{(m+1)}=\bigcup_{I\in\mathcal{T}_{(m)}%
}\mathcal{C}(I)$

(iii) The family $\mathcal{T}$ differentiates $L^{1}(X,\mu)$.

(b) Given an $N$-homogeneous tree-like family $\mathcal{T}$ on $X$ the
corresponding maximal operator $M_{\mathcal{T}}$ is defined for any $\psi\in
L^{1}(X,\mu)$ by%
\[
M_{\mathcal{T}}\psi(x)=\sup\left\{  \frac{1}{\mu(I)}\int_{I}\left\vert
\psi\right\vert d\mu:x\in I\in\mathcal{T}\right\}  \text{.}%
\]

\end{definition}

\bigskip The dyadic maximal operator localized to $X=[0,1]^{n}$ is contained
in the above definition, with $N=2^{n}$.

Sharp upper estimates for such (as well as in the much more general case where
the homogeneity of the tree is not assumed) operators have been provided by
the evaluation of corresponding Bellman functions in various cases (see
\cite{Burk1}, \cite{Burk2}, \cite{Mel1}, \cite{Mel2}, \cite{Mel3}, \cite{Sla},
\cite{Sla1}, \cite{Vas}, \cite{Vas1}, \cite{Vas2}, \cite{SlSt}).

In \cite{Mel1a} corresponding sharp lower $L^{p}$-bounds have been found for
this operator, when the $L^{1}$ and the $L^{p}$ norms of the function are
fixed, and this was done by proving that given any $p>1$ and positive real
numbers $f,F$ with $f^{p}\leq F$ we have%
\begin{align}
\inf\{\int_{X}(M_{\mathcal{T}}\phi)^{p}d\mu &  :\phi\geq0\text{ measurable,
}\int_{X}\phi d\mu=f\text{, \ }\int_{X}\phi^{p}d\mu=F\}=\nonumber\\
&  =f^{p}+\frac{N^{p}-1}{N^{p}-N}(F-f^{p}). \label{a1}%
\end{align}

The purpose of this paper is to study further lower bonds for these operators
where the $L^{p}$-norm size condition of $M_{\mathcal{T}}\phi$ is replaced by
other size conditions such as $L^{p}$-integral on subsets of fixed measure,
weak $L^{q}$ where $q>p$ type size conditions and strong $L^{q}$ with $q$
different from $p$, thus providing more information on the lower bounds for
such operators.

In this direction we first define the following Bellman type function%
\begin{gather}
\mathcal{D}_{p}^{\mathcal{T}}(F,f,\kappa)=\inf{\huge \{}\sup_{\mu(E)=\kappa
}\int_{E}(M_{\mathcal{T}}\phi)^{p}d\mu:\phi\geq0\text{ measurable
with}\nonumber\\
\int_{X}\phi d\mu=f\text{, \ }\int_{X}\phi^{p}d\mu=F{\huge \}}\text{.}%
\end{gather}
the inner supremum taken over all measurable subsets $E$ of $X$ having measure
$\kappa$, where $\kappa\in(0,1]$, and the positive numbers $F,f$ are such that
$f^{p}\leq F$. Then we will prove the following

\begin{theorem}
For any $N$-homogeneous tree-like family $\mathcal{T}$ any $p>1$ any $F,f$
with $f^{p}\leq F$ and any $\kappa\in(0,1]$ we have%
\begin{equation}
\mathcal{D}_{p}^{\mathcal{T}}(F,f,\kappa)=\min\{\kappa u^{p}+\frac{N^{p}%
-1}{N^{p}-N}(F-u^{p-1}f):f\leq u\leq\min((\frac{F}{f})^{1/(p-1)},\frac
{f}{\kappa})\}
\end{equation}
and writing $c(N,p)=\frac{p-1}{p}\frac{N^{p}-1}{N^{p}-N}<1$,
\begin{gather}
\mathcal{D}_{p}^{\mathcal{T}}(F,f,\kappa)=\label{e1}\\
=\left\{
\begin{array}
[c]{l}%
\kappa f^{p}+\frac{N^{p}-1}{N^{p}-N}(F-f^{p})\text{ if \ }c(N,p)\leq\kappa
\leq1\\
\frac{N^{p}-1}{N^{p}-N}(F-c(N,p)^{p-1}\frac{f^{p}}{p\kappa^{p-1}})\text{ if
\ }c(N,p)(\frac{f^{p}}{F})^{1/(p-1)}\leq\kappa\leq c(N,p)\\
\kappa(\frac{F}{f})^{p/(p-1)}\text{ if \ }0<\kappa\leq c(N,p)(\frac{f^{p}}%
{F})^{1/(p-1)}\text{.}%
\end{array}
\right.  \nonumber
\end{gather}

\end{theorem}

From the above theorem one obtains lower bounds for the following equivalent
norm on weak $L^{q}$ when $q>p$:
\[
\left\Vert \psi\right\Vert _{q,\infty}=\sup_{0<\mu(E)}\mu(E)^{-\frac{1}%
{p}+\frac{1}{q}}\left(  \int_{E}\left\vert \psi\right\vert ^{p}d\mu\right)
^{1/p}%
\]

\begin{corollary}
Given $q>p>1$ and $F,f>0$ with $f^{p}\leq F$ we have for any measurable
$\phi\geq0$ on $X$ with $\int_{X}\phi d\mu=f$, \ $\int_{X}\phi^{p}d\mu=F$ the following:

i) If $\frac{q-1}{q-p}\frac{f^{p}}{F}<1$
\[
\left\Vert M_{\mathcal{T}}\phi\right\Vert _{q,\infty}\geq\max[c(N,p)^{1/q}%
\frac{(q-p)^{(q-p)/q(p-1)}}{(q-1)^{(q-1)/q(p-1)}}(\frac{F^{\frac{q-1}{p-1}}%
}{f^{\frac{q-p}{p-1}}})^{1/q},(f^{p}+\frac{N^{p}-1}{N^{p}-N}(F-f^{p}%
))^{1/p}]\text{.}%
\]

ii) If  $\frac{q-1}{q-p}\frac{f^{p}}{F}\geq1$
\[
\left\Vert M_{\mathcal{T}}\phi\right\Vert _{q,\infty}\geq\max[c(N,p)^{1/q}%
(\frac{p}{p-1})^{1/p}(F-\frac{f^{p}}{p})^{1/p},(f^{p}+\frac{N^{p}-1}{N^{p}%
-N}(F-f^{p}))^{1/p}]\text{.}%
\]

\begin{proof}
Clearly $(\left\Vert M_{\mathcal{T}}\phi\right\Vert _{q,\infty})^{p}\geq
\sup\{\kappa^{-1+\frac{p}{q}}\mathcal{D}_{p}^{\mathcal{T}}(F,f,\kappa
):0<\kappa\leq1\}$. Computing the derivative of this function of $\kappa$ in
each of the ranges described in (\ref{e1}) it is easy to see that it is: a)
increasing in \ $0<\kappa\leq c(N,p)(\frac{f^{p}}{F})^{1/(p-1)}$, b) doesn't
have interior local maximum in $c(N,p)\leq\kappa\leq1$ and c) in
$c(N,p)(\frac{f^{p}}{F})^{1/(p-1)}\leq\kappa\leq c(N,p)$ it has an interior
local maximum at $\kappa_{0}=(\frac{q-1}{q-p}\frac{f^{p}}{F})^{1/(p-1)}c(N,p)$
if $\frac{q-1}{q-p}\frac{f^{p}}{F}<1$ (note that $\frac{q-1}{q-p}>1$) and is
increasing there otherwise (hence is maximized at $c(N,p)$). Thus the supremum
is attained either at $\kappa=\kappa_{0}$ or $c(N,p)$ (depending on
$\frac{q-1}{q-p}\frac{f^{p}}{F}$) or at $\kappa=1$. Introducing these values
in (\ref{e1}) completes the proof.
\end{proof}
\end{corollary}

Next we examine the strong $L^{q}$ norms considering the following Bellman
type function%
\[
\mathcal{B}_{p,q}^{\mathcal{T}}(F,f)=\inf\{\int_{X}(M_{\mathcal{T}}\phi
)^{q}d\mu:\phi\geq0\text{ measurable, }\int_{X}\phi d\mu=f\text{, \ }\int
_{X}\phi^{p}d\mu=F\}\text{.}%
\]
The case $p=q$ has been treated in \cite{Mel1a}. Here we prove first that:

\begin{proposition}
For any $N$-homogeneous tree-like family $\mathcal{T}$ any $p>1$ any $q<p$ and
any $F,f$ with $f^{p}\leq F$ we have%
\begin{equation}
\mathcal{B}_{p,q}^{\mathcal{T}}(F,f)=f^{q}\text{.}\label{a7}%
\end{equation}

\begin{proof}
It suffices to take a large integer $m$ and an $I\in\mathcal{T}_{(m)}$ (thus
of measure $N^{-m}$) choose a function $\eta$ on $I$ such that $\int_{I}\eta
d\mu=fN^{-m}$, $\int_{I}\eta^{p}d\mu=F-f^{p}(1-N^{-m})$ and $\int_{I}\eta
^{q}d\mu$ is sufficiently small (depending on $m$). For example one may take a
function of the form $\eta=a\chi_{C}$ for $C\subset I$ and $a>0$. Next taking
$\phi=\eta\chi_{I}+f\chi_{X\backslash I}$ we conclude that $\phi$ satisfies
$\int_{X}\phi d\mu=f$, \ $\int_{X}\phi^{p}d\mu=F$ and $M_{\mathcal{T}}\phi=f$
on $X\backslash I$ whereas $\int_{I}(M_{\mathcal{T}}\eta)^{q}d\mu\leq
c_{q}\int_{I}\eta^{q}$ will be small. Then $m\rightarrow\infty$ proves
(\ref{a7}).
\end{proof}
\end{proposition}

Thus the interesting case is when $q>p$ and in this case we will prove the following

\begin{theorem}
For any $N$-homogeneous tree-like family $\mathcal{T}$ any $p>1$ any $q>p$ and
any $F,f$ with $f^{p}\leq F$ we have%
\begin{equation}
\mathcal{B}_{p,q}^{\mathcal{T}}(F,f)\geq f^{q}+\frac{N^{q}-1}{N^{q}%
-N}{\Large (}\frac{F^{\frac{q-1}{p-1}}}{f^{\frac{q-p}{p-1}}}-f^{q}%
{\Large )}\label{a5}%
\end{equation}
and we have equality when $(F/f^{p})^{1/(p-1)}$ is a power of $N$, that is if
$m$ is a nonnegative integer then%
\begin{equation}
\mathcal{B}_{p,q}^{\mathcal{T}}(N^{m(p-1)}f^{p},f)=f^{q}{\Large [}%
1+\frac{N^{q}-1}{N^{q}-N}{\Large (}N^{m(q-1)}-1{\Large )]}\text{.}\label{a6}%
\end{equation}

\end{theorem}

In section 2 we prove Theorem 1 and in section 3 we prove Theorem 2.

\section{Proof of Theorem 1}

To prove Theorem 1 we fix a measurable $\phi\geq0$ with $\int_{X}\phi d\mu=f$,
\ $\int_{X}\phi^{p}d\mu=F$ and a $\kappa\in(0,1)$ and then find $u\geq f$ such
that%
\begin{equation}
\kappa_{1}=\mu(\{M_{\mathcal{T}}\phi>u\})\leq\kappa\leq\mu(\{M_{\mathcal{T}%
}\phi\geq u\})\leq\frac{f}{u}%
\end{equation}
(thus $u\leq\frac{f}{\kappa}$) and it is easy to see that%
\begin{equation}
\sup_{\mu(E)=\kappa}\int_{E}(M_{\mathcal{T}}\phi)^{p}d\mu=\int
_{\{M_{\mathcal{T}}\phi>u\}}(M_{\mathcal{T}}\phi)^{p}d\mu+(\kappa-\kappa
_{1})u^{p}\text{.}\label{a2}%
\end{equation}
Next we obviously have $\{M_{\mathcal{T}}\phi>u\}=%
{\displaystyle\bigcup\limits_{j}}
I_{j}$ for a certain family $\{I_{j}\}$ of pairwise disjoint elements of
$\mathcal{T}$ maximal under $\frac{1}{\mu(I_{j})}\int_{I_{j}}\phi>u$.

By writing%
\begin{equation}
\lambda_{j}=\mu(I_{j})\text{, }~\beta_{j}=\frac{1}{\mu(I_{j})}\int_{I_{j}}%
\phi\text{, }\alpha_{j}=\int_{I_{j}}\phi^{p}\text{\ }%
\end{equation}
considering the trees $\mathcal{T}(I_{j})=\{I\in\mathcal{T}:I\subset I_{j}\}$
on the probability spaces $(I_{j},\frac{1}{\lambda_{j}}\mu)$ and applying
(\ref{a1}) to them we get for each $j$%
\begin{equation}
\frac{1}{\lambda_{j}}\int_{I_{j}}(M_{\mathcal{T}}\phi)^{p}d\mu\geq\frac
{1}{\lambda_{j}}\int_{I_{j}}(M_{\mathcal{T}(I_{j})}\phi)^{p}d\mu\geq
\frac{N^{p}-1}{N^{p}-N}\frac{a_{j}}{\lambda_{j}}-\frac{N-1}{N^{p}-N}\beta
_{j}^{p}.
\end{equation}
Hence adding these inequalities we get with $A=%
{\textstyle\sum\limits_{j}}
\alpha_{j}$, $B=%
{\textstyle\sum\limits_{j}}
\lambda_{j}\beta_{j}$ and noting that $%
{\textstyle\sum\limits_{j}}
\lambda_{j}=\kappa_{1}$, the following%
\begin{gather}
\int_{\{M_{\mathcal{T}}\phi>u\}}(M_{\mathcal{T}}\phi)^{p}d\mu\geq\frac
{N^{p}-1}{N^{p}-N}%
{\textstyle\sum\limits_{j}}
\alpha_{j}-\frac{N-1}{N^{p}-N}%
{\textstyle\sum\limits_{j}}
\lambda_{j}\beta_{j}^{p}=\nonumber\\
=\kappa_{1}u^{p}+\frac{N^{p}-1}{N^{p}-N}(A-Bu^{p-1})-\nonumber\\
-%
{\textstyle\sum\limits_{j}}
\lambda_{j}(u^{p}-\frac{N^{p}-1}{N^{p}-N}\beta_{j}u^{p-1}+\frac{N-1}{N^{p}%
-N}\beta_{j}^{p})\geq\nonumber\\
\geq\kappa_{1}u^{p}+\frac{N^{p}-1}{N^{p}-N}(A-Bu^{p-1})
\end{gather}
the last inequality follows since the maximality of $I_{j}$'s imply that
$u<\beta_{j}\leq Nu$ and then the convexity of the function $h(t)=1-\frac
{N^{p}-1}{N^{p}-N}t+\frac{N-1}{N^{p}-N}t^{p}$ combined with the fact that
$h(1)=h(N)=0$ give $u^{p}-\frac{N^{p}-1}{N^{p}-N}\beta_{j}u^{p-1}+\frac
{N-1}{N^{p}-N}\beta_{j}^{p}\leq0$ for each $j$.

Next note that $A=\int_{\{M_{\mathcal{T}}\phi>u\}}\phi^{p}d\mu$,
$B=\int_{\{M_{\mathcal{T}}\phi>u\}}\phi d\mu$, so Holder's inequality gives
$B^{p}\leq\kappa_{1}^{p-1}A$ and also note that $\phi\leq M_{\mathcal{T}}%
\phi\leq u$ on $D=X\backslash\{M_{\mathcal{T}}\phi>u\}$ which gives combined
with
\begin{equation}
F-A=\int_{D}\phi^{p}\leq u^{p-1}\int_{D}\phi=u^{p-1}(f-B)\text{.}\label{b}%
\end{equation}
The inequalities $u<\beta_{j}\leq Nu$ on the other hand give that
$B=\kappa_{1}xu$ where $1<x\leq N$ and so $A\geq\kappa_{1}^{-p+1}B^{p}%
=\kappa_{1}x^{p}u^{p}$ which combined with (\ref{b}) gives%
\begin{equation}
A\geq\max(F-u^{p-1}(f-\kappa_{1}xu),\kappa_{1}x^{p}u^{p})
\end{equation}
and so since $x>1$%
\begin{gather}
\int_{\{M_{\mathcal{T}}\phi>u\}}(M_{\mathcal{T}}\phi)^{p}d\mu\geq\nonumber\\
\geq\kappa_{1}u^{p}+\frac{N^{p}-1}{N^{p}-N}(\max(F-u^{p-1}(f-\kappa
_{1}xu),\kappa_{1}x^{p}u^{p})-\kappa_{1}xu^{p})=\nonumber\\
=\kappa_{1}u^{p}+\frac{N^{p}-1}{N^{p}-N}\max(F-u^{p-1}f,\kappa_{1}%
(x^{p}-x)u^{p})\geq\nonumber\\
\geq\kappa_{1}u^{p}+\frac{N^{p}-1}{N^{p}-N}\max(F-u^{p-1}f,0)\text{.}%
\end{gather}
Now this combined with (\ref{a2}) gives%
\begin{gather}
\sup_{\mu(E)=\kappa}\int_{E}(M_{\mathcal{T}}\phi)^{p}d\mu\geq\kappa
u^{p}+\frac{N^{p}-1}{N^{p}-N}\max(F-u^{p-1}f,0)\geq\nonumber\\
\geq\min\{\kappa u^{p}+\frac{N^{p}-1}{N^{p}-N}(F-u^{p-1}f):f\leq u\leq
(\frac{F}{f})^{1/(p-1)}\text{ and }u\leq\frac{f}{\kappa}\}\text{.}\label{a3}%
\end{gather}

Conversely given $F,f,\kappa$ as above we let $u_{0}\in\lbrack f,\min
((\frac{F}{f})^{1/(p-1)},\frac{f}{\kappa})]$ minimize the quantity $\kappa
u^{p}+\frac{N^{p}-1}{N^{p}-N}(F-u^{p-1}f)$ in the above inequality and using
Lemma 1 and the proof of Proposition 1 in \cite{Mel1a} we can find pairwise
disjoint elements $\{I_{j}\}$ of $\mathcal{T}$ and measurable functions
$\phi_{j}\geq0$ on each $I_{j}$ such that%
\[
\sum_{j}\mu(I_{j})=\kappa\text{, }\frac{1}{\mu(I_{j})}\int_{I_{j}}\phi
_{j}=u_{0}\text{, }\int_{I_{j}}\phi_{j}^{p}=\alpha_{j}\geq\mu(I_{j})u_{0}%
^{p}\text{ }%
\]
with%
\[
\sum_{j}\alpha_{j}=F-u_{0}^{p-1}f+\kappa u_{0}^{p-1}\geq\kappa u_{0}^{p-1}%
\]
and such that for each $j$%
\[
\frac{1}{\mu(I_{j})}\int_{I_{j}}(M_{\mathcal{T}(I_{j})}\phi_{j})^{p}d\mu
=u_{0}^{p}+\frac{N^{p}-1}{N^{p}-N}(\frac{\alpha_{j}}{\mu(I_{j})}-u_{0}%
^{p})\text{.}%
\]
Now let $Y=X\backslash%
{\textstyle\bigcup\limits_{j}}
I_{j}$ and choose a measurable $P\subset Y$ such that $\mu(P)=\dfrac{f-\kappa
u_{0}}{u_{0}}\in\lbrack0,1-\kappa]$ by the conditions in (\ref{a3}) and define
the measurable function%
\begin{equation}
\phi=u_{0}\chi_{P}+\sum_{j}\phi_{j}\chi_{I_{j}}\text{.}%
\end{equation}
Since $\int_{Y}\phi=f-\kappa u_{0}$ and $\int_{Y}\phi^{p}=u_{0}^{p-1}f-\kappa
u_{0}^{p}=F-\sum_{j}\alpha_{j}$ we easily obtain that $\int_{X}\phi d\mu=f$,
\ $\int_{X}\phi^{p}d\mu=F$. However since $\phi\leq u_{0}$ on $Y$ and
$\frac{1}{\mu(I_{j})}\int_{I_{j}}\phi=u_{0}$ we conclude that $M_{\mathcal{T}%
}\phi=M_{\mathcal{T}(I_{j})}\phi_{j}$ on each $I_{j}$ and that
$\{M_{\mathcal{T}}\phi>u_{0}\}\subset%
{\textstyle\bigcup\limits_{j}}
I_{j}\subset\{M_{\mathcal{T}}\phi\geq u_{0}\}$ hence$^{{}}$%
\begin{gather}
\sup_{\mu(E)=\kappa}\int_{E}(M_{\mathcal{T}}\phi)^{p}d\mu=\int_{%
{\textstyle\bigcup\limits_{j}}
I_{j}}(M_{\mathcal{T}}\phi)^{p}d\mu=\nonumber\\
=\sum_{j}\mu(I_{j})[u_{0}^{p}+\frac{N^{p}-1}{N^{p}-N}(\frac{\alpha_{j}}%
{\mu(I_{j})}-u_{0}^{p})]=\nonumber\\
=\kappa u_{0}^{p}+\frac{N^{p}-1}{N^{p}-N}(F-u_{0}^{p-1}f)\text{.}%
\end{gather}
To prove (\ref{e1}) we first note that $c(N,p)<1$ since $N^{p}-1>p(N-1)$. Next
by defining $g(u)=\kappa u^{p}+\frac{N^{p}-1}{N^{p}-N}(F-u^{p-1}f)$ we observe
that $g^{\prime}(u)=0$ iff $u=c(N,p)\frac{f}{\kappa}$ and that $g(t^{1/(p-1)}%
)$ is convex on $t>0$. Hence to find the minimum of $g(u)$ in $[f,\min
((\frac{F}{f})^{1/(p-1)},\frac{f}{\kappa})]$ it suffices to examine the
relative position of the values $f,(\frac{F}{f})^{1/(p-1)},\frac{f}{\kappa}$
and $~c(N,p)\frac{f}{\kappa}$. Thus when $c(N,p)\frac{f}{\kappa}\leq f$ that
is $\kappa\geq c(N,p)$ the minimum is attained for $u=f$, when $f<c(N,p)\frac
{f}{\kappa}\leq(\frac{F}{f})^{1/(p-1)}$ that is $c(N,p)(\frac{f^{p}}%
{F})^{1/(p-1)}\leq\kappa\leq c(N,p)$ the minimum is attained at $u=c(N,p)\frac
{f}{\kappa}$ and when $0<\kappa<c(N,p)(\frac{f^{p}}{F})^{1/(p-1)}$ the minimum
is attained for $u=(\frac{F}{f})^{1/(p-1)}$. Substituting the corresponding
values of $u$ in $g$ we obtain (\ref{e1}). This completes the proof of Theorem 1.

\section{Proof of Theorem 2}

Here we will prove Theorem 2 where $q>p>1$. For the lower bound we follow a
classical Bellman type argument. Assuming that $\mathcal{T}$ is a
$N$-homogeneous tree define the following function %

\begin{gather}
\mathcal{B}_{p,q}^{\mathcal{T}}(F,f,L)=\sup{\huge \{}\int_{X}\max
(M_{\mathcal{T}}\phi,L)^{q}d\mu:\phi\geq0\text{ measurable with}\nonumber\\
\int_{X}\phi d\mu=f\text{, \ }\int_{X}\phi^{p}d\mu=F{\huge \}}\text{.}%
\end{gather}
whenever $F,f,L$ are positive real numbers with $f\leq L$ and $f^{p}\leq F$.

Then we will prove the following from which (\ref{a5}) easily follows (by
taking $L=f$).

\begin{lemma}
We have
\begin{equation}
\mathcal{B}_{p,q}^{\mathcal{T}}(F,f,L)\geq L^{q}+\frac{N^{q}-1}{N^{q}-N}%
(\frac{F^{\frac{q-1}{p-1}}}{f^{\frac{q-p}{p-1}}}-L^{q-1}f)^{+} \label{b1}%
\end{equation}
where $x^{+}=\max(x,0)$.

\begin{proof}
Write $r=\dfrac{q-1}{p-1}>1$. We first consider a nonnegative $\mathcal{T}%
$-step function at level $m$, $\phi\geq0$ that is a finite linear combination
of the functions $\chi_{I}$ where $I\in\mathcal{T}_{(m)}$, such that $\int
_{X}\phi d\mu=f$ and $\int_{X}\phi^{p}d\mu=F$ and prove (\ref{b1}) by
induction on $m$, the case $m=0$ being trivial.

We have $X=%
{\textstyle\bigcup\limits_{i=1}^{N}}
J_{i}$ where each $J_{i}$ is in $\mathcal{T}_{(1)}$ and has measure $1/N$ and
we write%
\[
F_{i}=\frac{1}{\mu(J_{i})}\int_{J_{i}}\phi^{p}\text{, }f_{i}=\frac{1}%
{\mu(J_{i})}\int_{J_{i}}\phi\text{.}%
\]
Note that the restriction of $\phi$ on each $J_{i}$ is a $\mathcal{T}(J_{i}%
)$-step function at level $m-1$. Also in the case $f_{i}>L$ we have
$\max(M_{\mathcal{T}}\phi,L)=M_{\mathcal{T}(J_{i})}\phi$ on $J_{i}$ and in the
case $f_{i}\leq L$ we have $\max(M_{\mathcal{T}}\phi,L)=\max(M_{\mathcal{T}%
(J_{i})}\phi,L)$ on $J_{i}$. hence by the induction hypothesis we have%
\begin{gather*}
N\int_{X}\max(M_{\mathcal{T}}\phi,L)^{q}d\mu=\\
=\left(  \sum_{f_{i}\leq L}\int_{X}\max(M_{\mathcal{T}(J_{i})}\phi,L)^{q}%
\frac{d\mu}{\mu(J_{i})}+\sum_{f_{i}>L}\int_{X}(M_{\mathcal{T}(J_{i})}\phi
)^{q}\frac{d\mu}{\mu(J_{i})}\right)  \geq\\
\geq\sum_{f_{i}\leq L}(L^{q}+\frac{N^{q}-1}{N^{q}-N}(\frac{F_{i}^{r}}%
{f_{i}^{r-1}}-L^{q-1}f_{i})^{+})+\sum_{f_{i}>L}(f_{i}^{q}+\frac{N^{q}-1}%
{N^{q}-N}(\frac{F_{i}^{r}}{f_{i}^{r-1}}-f_{i}^{q}))\text{.}%
\end{gather*}
Next we observe that when $f_{i}>L$ we also have $f_{i}\leq Nf\leq NL$, thus
by the convexity of the function $h(t)=1-\frac{N^{p}-1}{N^{p}-N}t+\frac
{N-1}{N^{p}-N}t^{p}$ since $h(1)=h(N)=0$ and since $L<f_{i}\leq(F_{i}%
/f_{i})^{1/(p-1)}~$we have
\[
f_{i}^{q}+\frac{N^{q}-1}{N^{q}-N}(\frac{F_{i}^{r}}{f_{i}^{r-1}}-f_{i}^{q})\geq
L^{q}+\frac{N^{q}-1}{N^{q}-N}(\frac{F_{i}^{r}}{f_{i}^{r-1}}-L^{q-1}f_{i}%
)^{+}\text{.}%
\]
Therefore using the inequality $(a_{1}+...+a_{N})^{+}\leq a_{1}^{+}%
+...+a_{N}^{+}$ we get
\begin{align*}
\int_{X}\max(M_{\mathcal{T}}\phi,L)^{q}d\mu &  \geq\frac{1}{N}\sum_{i=1}%
^{N}(L^{q}+\frac{N^{q}-1}{N^{q}-N}(\frac{F_{i}^{r}}{f_{i}^{r-1}}-L^{q-1}%
f_{i})^{+})\geq\\
&  \geq L^{q}+\frac{N^{q}-1}{N^{q}-N}(\frac{1}{N}\sum_{i=1}^{N}\frac{F_{i}%
^{r}}{f_{i}^{r-1}}-L^{q-1}f)^{+}%
\end{align*}
since $Nf=f_{1}+...+f_{N}$. Now using Holder's inequality for $r>1$ we have
\[
\frac{1}{N}\sum_{i=1}^{N}\frac{F_{i}^{r}}{f_{i}^{r-1}}\geq\frac{1}{N}%
\frac{(\sum_{i=1}^{N}F_{i})^{r}}{(\sum_{i=1}^{N}f_{i})^{r-1}}=\frac{F^{r}%
}{f^{r-1}}%
\]
and this completes the induction.

For the general case, given $\phi\geq0$ measurable satisfying $\int_{X}\phi
d\mu=f$ and $\int_{X}\phi^{p}d\mu=F$ \ we define $\phi_{m}$ as follows%
\[
\phi_{m}=\sum_{I\in\mathcal{T}_{(m)}}\operatorname{Av}_{I}(\phi)\chi_{I}%
\]
and\ we note that
\begin{equation}
M_{\mathcal{T}}\phi_{m}=\sum\limits_{I\in\mathcal{T}_{(m)}}\max
\{\operatorname{Av}_{J}(\phi):I\subseteq J\in\mathcal{T}\}\chi_{I}%
\end{equation}
since $\operatorname{Av}_{J}(\phi)=\operatorname{Av}_{J}(\phi_{m})$ whenever
$I\subseteq J\in\mathcal{T}$ , $I\in\mathcal{T}_{(m)}$. Also
\begin{equation}
\int_{X}\phi_{m}d\mu=\int_{X}\phi d\mu=f\text{, \ }F_{m}=\int_{X}\phi_{m}%
^{p}d\mu\leq\int_{X}\phi^{p}d\mu=F
\end{equation}
for all $m$ and $M_{\mathcal{T}}\phi_{m}$ converges monotonically to
$M_{\mathcal{T}}\phi$. Also since each $\phi_{m}$ is a $\mathcal{T}$-step
function for every $m$ the inequality (\ref{b1}) holds for each $\phi_{m}$
with $F_{m}$ in the place of $F$. On the other hand we have $\phi_{m}^{p}%
\leq(M_{\mathcal{T}}\phi)^{p}$ everywhere and $\phi_{m}^{p}\rightarrow\phi
^{p}$ almost everywhere by property (iv) in Definition 1. Hence by dominated
convergence we conclude that $F_{m}=\int_{X}\phi_{m}^{p}d\mu\rightarrow
\int_{X}\phi^{p}d\mu=F$ and so using monotone convergence for $M_{\mathcal{T}%
}\phi_{m}$ we easily get (\ref{b1}) for $\phi$.
\end{proof}
\end{lemma}

Now in the case where $f=1$ and $F=N^{m(p-1)}$ we let $X=I_{0}\supseteq
I_{1}\supseteq...I_{s}\supseteq I_{s+1}\supseteq...\supseteq I_{m}$ be a chain
such that $I_{s}\in\mathcal{T}_{(s)}$ for all $s$ (and so $\mu(I_{s})=N^{-s}$)
and consider the function%
\begin{equation}
\phi=N^{m}\chi_{I_{m}}\label{c1}%
\end{equation}
which clearly satisfies $\int_{X}\phi d\mu=f$ and $\int_{X}\phi^{p}d\mu=F$ and
as it is easy to see that%
\begin{equation}
M_{\mathcal{T}}\phi=N^{m}\chi_{I_{m}}+N^{m-1}\chi_{I_{m-1}\backslash I_{m}%
}+....+N\chi_{I_{1}|I_{0}}+\chi_{I_{0}|I_{1}}\label{c2}%
\end{equation}
we get with $r=\dfrac{q-1}{p-1}$
\begin{align}
\int_{X}(M_{\mathcal{T}}\phi)^{q}d\mu &  =N^{(q-1)m}+(1-\frac{1}%
{N})(N^{(q-1)(m-1)}+...+N^{q-1}+1)=\nonumber\\
&  =N^{(q-1)m}+(1-\frac{1}{N})\frac{N^{(q-1)m}-1}{N^{q-1}-1}=\nonumber\\
&  =1+\frac{N^{q}-1}{N^{q}-N}(F^{r}-1)\text{.}\label{c3}%
\end{align}
Therefore by homogeneity we conclude that (\ref{a5}) is an equality when
$(F/f^{p})^{1/(p-1)}$ is a power of $N$ and this proves (\ref{a6}).


\begin{thebibliography}{99}                                                                                               %


\bibitem {Burk1}D. L. Burkholder, Martingales and Fourier analysis in Banach
spaces, C.I.M.E. Lectures (Varenna (Como), Italy, 1985), \textit{Lecture Notes
in Mathematics} 1206 (1986), 61-108.

\bibitem {Burk2}D. L. Burkholder, Boundary value problems and sharp
inequalities for martingale transforms, \textit{Ann. of Prob.} \textbf{12}
(1984), 647-702.

\bibitem {Mel1}A. D. Melas, The Bellman functions of dyadic-like maximal
operators and related inequalities, \textit{Adv. in Math. }\textbf{192}
(2005), 310-340.

\bibitem {Mel1a}A. D. Melas, E. Nikolidakis, Th. Stavropoulos, Sharp local
$\ $lower $L^{p}$-bounds for dyadic-like maximal operators, \textit{Proc..
Amer. Math. Soc. }\textbf{141} (2013), 3171-3181.

\bibitem {Mel2}A. D. Melas, Sharp general local estimates for Dyadic-like
maximal operators and related Bellman functions, \textit{Adv. in Math.
}\textbf{220}, no. 2, (2009), 367-426.

\bibitem {Mel3}A. D. Melas, Dyadic-like maximal operators on $L\log L$
functions, \textit{Jour. Funct. Anal}. \textbf{257}, no. 6, (2009), 1631-1654.

\bibitem {Naz}F. Nazarov, S. Treil, The hunt for a Bellman function:
applications to estimates for singular integral operators and to other
classical problems of harmonic analysis, \textit{Algebra i Analyz }\textbf{8}
no. 5 (1996), 32-162.

\bibitem {Naz1}F. Nazarov, S. Treil, A. Volberg, The Bellman functions and
two-weight inequalities for Haar multipliers. \textit{Journ. Amer}.
\textit{Math. Soc.} \textbf{12 }no. 4 (1999), 909-928.

\bibitem {Naz2}F. Nazarov, S. Treil, A. Volberg, Bellman function in
Stochastic Optimal Control and Harmonic Analysis (how our Bellman function got
its name), \textit{Oper. Theory: Advances and Appl. }\textbf{129} (2001),
393-424, Birkh\"{a}user Verlag.

\bibitem {SlSt}L. Slavin, A. Stokolos, V. Vasyunin. Monge-Amp\`{e}re equations
and Bellman functions: The dyadic maximal operator. \textit{C. R. Acad. Paris
S\'{e}r. I Math.}, \textbf{346}, no 9-10, (2008), 585-588.

\bibitem {Sla}L. Slavin, V. Vasyunin, Sharp results in the integral-form
John-Nirenberg inequality, \textit{Trans. Amer. Math. Soc.} \textbf{363}, no
8, (2011), 4135-4169.

\bibitem {Sla1}L. Slavin, A. Volberg. The explicit BF for a dyadic
Chang-Wilson-Wolff theorem. The s-function and the exponential integral.
\textit{Contmp. Math.,} \textbf{444} (2007), 215-228.

\bibitem {Ste}E. M. Stein. Note on the class $L\log L$, \textit{Studia Math}.
\textbf{32 }(1969), 305-310.

\bibitem {Vas}V. Vasyunin. The sharp constant in the reverse Holder inequality
for Muckenhoupt weights. \textit{Algebra i Analiz}, \textbf{15} (2003), no. 1, 73-117

\bibitem {Vas1}V. Vasyunin, A. Volberg. The Bellman functions for a certain
two weight inequality: the case study. \textit{Algebra i Analiz}, \textbf{18}
(2006), No. 2

\bibitem {Vas2}V. Vasyunin, A. Volberg. Monge-Ampere equation and Bellman
optimization of Carleson embedding Theorem, \textit{Amer. Math. Soc.
Translations Series }2, \textbf{226}, (2009).

\bibitem {Vol}A. Volberg, Bellman approach to some problems in Harmonic
Analysis, \textit{Seminaire des Equations aux deriv\'{e}es partielles, Ecole
Polyt\'{e}chnique, }2002, epos\'{e}. XX, 1-14.\textit{\ }

\bibitem {Wang}G. Wang, Sharp maximal inequalities for conditionally symmetric
martingales and Brownian motion, \textit{Proc. Amer. Math. Soc. }\textbf{112
}(1991), 579-586.
\end{thebibliography}
\end{document}